\documentclass{article}

\usepackage{latexsym}
\usepackage{amssymb}
\usepackage{amsmath,amsfonts}
\usepackage{graphicx}
\usepackage{tikz}
\usetikzlibrary{arrows}
\usepackage{enumerate}

\newtheorem{theorem}{Theorem}[section]
\newtheorem{prop}[theorem]{Proposition}
\newtheorem{lemma}{Lemma}[section]
\newtheorem{corollary}{Corollary}[section]
\newtheorem{example}{Example}[section]
\newtheorem{definition}{Definition}[section]
\newtheorem{remark}{Remark}[section]
\newtheorem{conjecture}[theorem]{Conjecture}

\newtheorem{question}[theorem]{Question}

\def \dd{\hfill \baseb \vskip .5cm}
\def \d{{\noindent \it Proof. } }

\def\ex{\begin{example}}
\def\eex{\end{example}}
\def\exx{\end{example}}
\def\t{\begin{theorem}}
\def\tt{\end{theorem}}
\def\D{\begin{definition}}
\def\DD{\end{definition}}
\def\l{\begin{lemma}}
\def\ll{\end{lemma}}
\def\c{\begin{corollary}}
\def\cc{\end{corollary}}
\def\cj{\begin{conjecture}}
\def\cjj{\end{conjecture}}
\def\e{\begin{equation}}
\def\ee{\end{equation}}
\def\p{\begin{prop}}
\def\pp{\end{prop}}
\def\q{\begin{question}}
\def\qq{\end{question}}

\def \Gn{{\mathcal G}_n}

\def \X{\mathcal X}
\def \Z{\mathcal Z}
\def\A{\mathcal A}
\def\B{\mathcal B}
\def\dom{\sqsubseteq}
\def\doms{\sqsupseteq}
\newcommand{\baseb}{\hfill \rule{2mm}{2mm}}

\begin{document}

\baselineskip 15pt


\title{ Graph Cordiality -- Extremes and Preservers\thanks{
\textbf{AMS Classification numbers}: 05C20, 05C60, 05C78} }

\author{LeRoy B. Beasley\\  \\Clocktower Plaza\#317 \\  550 North Main, Box C3 \\  Logan, Utah 84321, USA,\\
\text{leroy\_beas@aol.com}}




\maketitle

\begin{abstract}
An undirected graph is said to be cordial if there is a friendly (0,1)-labeling of the vertices that induces a friendly (0,1)-labeling of the edges. An undirected  graph $G$  is said to be $(2,3)$-orientable if there exists a friendly (0,1)-labeling of the vertices of $G$ such that about one third of the edges are incident to vertices labeled the same.  That is,  there is some digraph that is an  orientation of $G$ that  is $(2,3)$-cordial.  Examples of the smallest noncordial/non-$(2,3)$-orientable  graphs are given and upper bounds on the possible number of edges in a cordial/$(2,3)$-orientable graph are presented. It is also shown that if  $T$ is a linear operator on the set of all undirected graphs on $n$ vertices that strongly preserves the set of cordial graphs or the set of $(2,3)$-orientable graphs then $T$ is a vertex permutation..\

\end{abstract}

\section{Introduction}

Let $\Gn$ denote the set of all simple loopless undirected  graphs on the vertex set $V=\{v_1,v_2,\dots, v_n\}$.   Let $G=(V,E)$ be an undirected graph in $\Gn$ with vertex set $V$ and edge set $E$.  A $(0,1)$-labeling of the vertex set is a mapping $f:V\to \{0,1\}$ and is said to be {\em friendly} if approximately one half of the vertices are labeled 0 and the others labeled 1.  An induced labeling of the edge set  is a mapping $g:E\to \{0,1\}$ where for an edge $uv, g(uv)= \hat{g}(f(u),f(v)$ for some $\hat{g}:\{0,1\}\times\{0,1\}\to \{0,1\}$ and is said to be cordial if $f$ is friendly and about one half the edges of $G$ are labeled 0 and the others labeled 1, that is, both $f$ and $g$ are friendly.  A graph, $G$,  is called {\em cordial} if there exists a cordial induced labeling of the edge set of $G$. 

Let $|{\cal X}|$ denote the cardinality of the set ${\cal X}$, so that $|g^{-1}(0)|$ denotes the number of edges labeled $0$ if $g:E\to \{0,1\}$ is the labeling of the edges.  

 In this article we consider not only cordial undirected graphs but also use a cordial labeling of directed graphs that is not  merely a cordial labeling of the underlying undirected graph.  It was defined and studied in \cite{B}. 

We shall consider three classes of graphs:  Those that are cordial with an induced edge labeling consisting of the sum modulo 2 of the labels of the end vertices (see \cite{C});  those that are cordial with an induced edge labeling consisting of the product  of the labeling of the end vertices (see \cite{S}) and those that are $(2,3)$-orientable (see\cite{B}), that is, graphs that may be oriented such that if each arc  of the directed graph is labeled as the   label  of the terminal vertex minus the label of the initial vertex, and  we have that  about one third of the arcs are labeled 0, about one third labeled 1 and the remaining one third labeled -1.

\section{Preliminaries.}

We begin with formal definitions of the terms from the introduction.

\D (See \cite{H}) Let $\Z$ be a set and label the set with entries from the set  ${\cal A}$. Let $f$  the function $f:\Z\to {\cal A}$ that preforms this labeling.  This labeling of a set is called {\em\underline{$\A$-friendly}} if $-1\leq |f^{-1}(i)| - |f^{-1}(j)|  \leq 1$ for any $i,j\in{\A}$.  A labeling $f$ is said to be {\em\underline{$k$-friendly}} if it is $\Z_k$-friendly where $\Z_k$ is the set of integers modulo $k$.  Further the labeling $f$  is said to be {\em\underline{friendly}} if it is $2$-friendly. \DD

\D  Let $\A$ be a semiring and  $G$ be an undirected graph with vertex set $V$ and edge set $E$.  Let $f:V\to\A$ be an $\A$-friendly labeling of the \underline{non-isolated} vertices of $V$ and let $g:E\to\A$ be an induced labeling of $E$, that is given an edge $uv$, $g(uv) = \hat{g}(f(u),f(v))$  where $\hat{g}:\A\times\A\to \A$.  Note that $\hat{g}$ must be symmetric to be well defined.  
 The graph $G$ with the $\A$-friendly  labeling $f$ of the non-isolated vertices of $G$ and induced edge labeling $g$ is said to be {\em\underline{$\A$-cordial }} if $g$ is also an $\A$-friendly mapping.  Also, $G$ is {\em\underline{$k$-cordial}} if $\A$ is $\Z_k$, and we say $G$ is {\em\underline{cordial}} if $\A$ is $\Z_2$.
\DD
\begin{remark}   It should be noted that $(\Z_2,+)$-cordial and $(\Z_2,-)$-cordial are the same and in fact is the usual definition of cordial where $f$ is a  $(0,1)$ labeling and $g(u,v)=|f(u)-f(v)|$.   Further, the restriction that $f$ is a friendly labeling of the non-isolated vertices, not necessarily all vertices, is required because, for example, if the restriction to  non-isolated vertices were not required, the graph $2K_2$ would be   $(\Z_2,+)$-cordial as a graph on five vertices and not on four vertices. 
\end{remark}

For digraphs we have a similar definition, however, as digraphs are not usually symmetric, we do not require the induced mapping $g$ to be symmetric:
\D   Let $\A$ and $\B$ be sets and let $G=(V,A)$ be an directed graph with vertex set $V$ and arc set $A$.  Let $f:V\to\A$ be an $\A$-friendly labeling of the non-isolated vertices of  $V$ and let $g:A\to\B$ be an induced labeling of $A$, that is given an arc $\overrightarrow{uv}$, $g(uv) = \hat{g}(f(u),f(v))$  where $\hat{g}:\A\times\A\to \B$. 
 The graph $G$ with the $\A$-friendly  labeling $f$ of the non-isolated vertices  and induced edge labeling $g$ is said to be {\em\underline{$(\A,\B)$-cordial }} if $g$ is  a $\B$-friendly mapping. When $\A=\Z_k$ and $\B=\Z_\ell$ we say that $G$ is $(k,\ell)$-cordial.  In particular we say that $G$ is $(2,3)$-cordial if $\A=\Z_2$, $\B=\Z_3$ and $g(\overrightarrow{uv}) =f(v)-f(u)$.  Note that in this case $g$ is anti-symmetric.\DD

Note that in the above definitions, if $\hat{g}$ is the binary  mapping corresponding to one of the binary operations on $\A$ or $\B$ we indicate that by placing that operator in the notation.  For example a digraph is $(\A,\B,-)$-cordial indicates that the induced labeling $g$ is $g(\overrightarrow{uv})=f(v)-f(u)$.  So a digraph is $(2,3,-)$-cordial means that the arc labelings are $f(v)-f(u)$ for the arc $\overrightarrow{uv}$.  In this specific case we usually drop the minus sign and write $(2,3)$-cordial.  A graph with an $\A$-friendly vertex labeling is {\em\underline{product cordial}} if it is $(\A,\times)$-cordial (see \cite{S}).  Unless specified otherwise, product cordial graphs have the vertex set labeled with $\{0,1\}$.

\D An undirected graph $G$ is said to be {\em\underline{$(2,3)$-orientable}} if some orientation of the edges yeilds a $(2,3)$-cordial digraph.\DD

\newpage
\section{Extremes of Cordial Graphs}
In this section we shall exhibit examples of graphs that are the smallest or largest (edgewise) of cordial/noncordial graphs.
\subsection{ $(Z_2,+)$-Cordial/Non-cordial graphs} \label{Z+C}
\subsubsection{Smallest $(Z_2,+)$-Cordial/Non-cordial Graphs}  It is easily shown  that the graph with one edge and all graphs with three edges are  $(Z_2,+)$-cordial.  However, the connected graph with two edges (a 2-path or a 2-star)  is  $(Z_2,+)$-cordial while the graph $2K_2$ is not and the complete graph on $n\geq 4$ vertices  is not  $(Z_2,+)$-cordial..

\subsubsection{Largest $(Z_2,+)$-Cordial/Non-cordial Graphs}\label{L+}  If $n\geq 4$, $K_n$ is not  $(Z_2,+)$-cordial.  So how large can a graph be and be  $(Z_2,+)$-cordial?  That is:  How many edges can a $(Z_2,+)$-cordial graph on $n$ vertices  contain?  

Let $f$ be a friendly labeling of $K_n$.  Since the labeling of the vertices of a complete graph is independent of the visual location of the vertex, we may assume that the vertices labeled 0 are located on the right half and the vertices labeled 1 are on the left.  The vertices labeled 0 induce a graph of order $\frac{\lfloor \frac{n}{2}\rfloor(\lfloor \frac{n}{2}\rfloor-1)}{2}$  and the vertices labeled 1 induce a subgraph of order $\frac{\lceil \frac{n}{2}\rceil(\lceil \frac{n}{2}\rceil-1)}{2}$, or visa versa.  In either case there are precisely  $\frac{\lfloor \frac{n}{2}\rfloor(\lfloor \frac{n}{2}\rfloor-1)}{2}+\frac{\lceil \frac{n}{2}\rceil(\lceil \frac{n}{2}\rceil-1)}{2}$ edges labeled 0 and the others  labeled 1.  

So if $n$ is  even, say $n=2k$,  there are $k^2-k$ edges labeled 0 and $k^2$ edges labeled 1. Thus to have a $(Z_2,+)$-cordial graph one must delete at least $k-1$ edges joining vertices that are labeled differently.  Thus, when $n=2k$ is even,  there can be at most $2k^2-2k+1$ edges in a   $(Z_2,+)$-cordial graph.  

If $n$ is odd, say $n=2k+1$, there are $k^2+k$ edges labeled 1 and $k^2$ edges labeled 0.  Thus to get a graph that is  $(Z_2,+)$-cordial one must delete at least $k-1$ edges connecting vertices labeled differently.  Thus, when $n=2k+1$ is odd,  there can be at most $2k^2+1$ edges in a  $(Z_2,+)$-cordial graph.  

The added ``+1'' in the formulas accounts for the fact that for a cordial graph the cardinality of the two sets $g^{-1}(0)$ and $g^{-1}(1)$  may differ by 1.

\subsection{Product Cordial/Non-cordial Graphs.}
\subsubsection{Smallest $(\Z_2,\times)$-Cordial/Non-cordial Graphs}  It is easily shown  that all graphs with one, two, or three    edges are  $(Z_2,\times)$-cordial.  There are two graphs with four edges that are not  $(Z_2,\times)$-cordial, the four cycle and a three cycle with an attached edge.  .  
\subsubsection{Largest  $(Z_2,\times)$-cordial/Non-cordial Graphs.}  For a labeling of the vertices of a graph on $n$ vertices to be friendly, at most $\lceil\frac{n}{2}\rceil$ vertices can be labeled 1.  Thus, at most $\frac{\lceil\frac{n}{2}\rceil(\lceil\frac{n}{2}\rceil-1)}{2}$  edges can be labeled 1 for any product cordial graph on $n$ vertices.  Thus, a   $(\Z_2,\times)$-cordial graph on $n$ vertices can have at most at most $\lceil\frac{n}{2}\rceil(\lceil\frac{n}{2}\rceil-1)$ edges since $\frac{\lceil\frac{n}{2}\rceil(\lceil\frac{n}{2}\rceil-1)}{2}<\frac{1}{2}\frac{n(n-1)}{2}-1$  if $n\geq 4$.

\subsection{$(2,3)$-Orientable/Non-orientable Graphs.}

\subsubsection{Smallest $(2,3)$-Orientable/Non-orientable Graphs}  

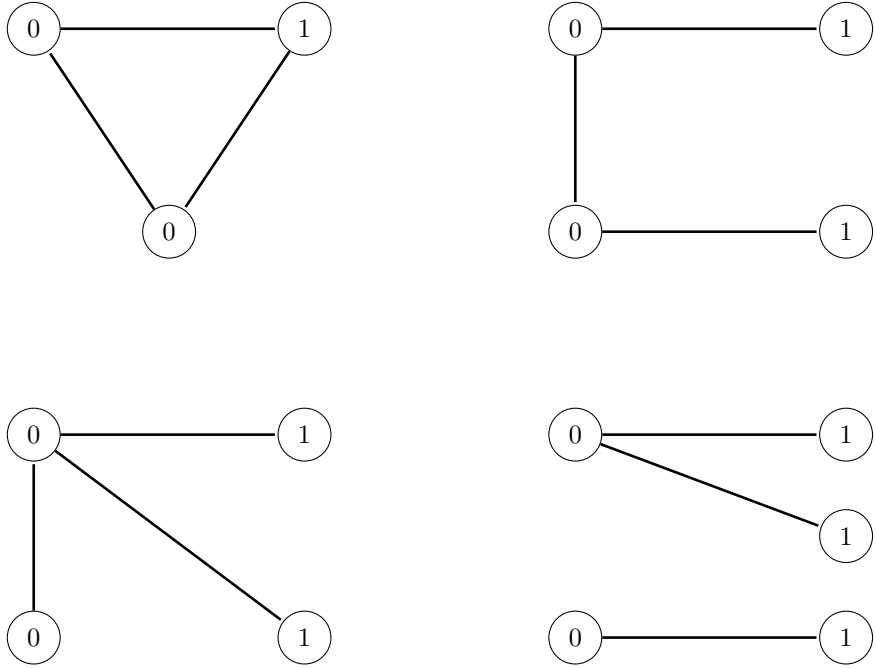
\begin{figure}[h]
\begin{center}
\begin{tikzpicture}[scale=0.9]

\tikzset{vertex/.style = {shape=circle,draw,minimum size=2em}}
\tikzset{edge/.style = {->,> = stealth',shorten >=1pt,thick}}

\node[vertex] (v1) at  (0,6) {$0$};
\node[vertex] (v2) at  (2,3) {$0$};
\node[vertex] (v3) at  (4,6) {$1$};
\node[vertex] (v2k-2) at  (8,6) {$0$};
\node[vertex] (v2k-1) at  (8,3) {$0$};
\node[vertex] (v2k) at  (12,6) {$1$};
\node[vertex] (v2k-3) at  (12,3) {$1$};

\draw[edge,-, line width=1.0pt] (v1) to (v3);
\draw[edge,-, line width=1.0pt] (v3) to (v2);
\draw[edge,-, line width=1.0pt] (v2) to (v1);
\draw[edge,-, line width=1.0pt] (v2k-2) to (v2k-1);
\draw[edge,-, line width=1.0pt] (v2k-2) to (v2k);
\draw[edge,-, line width=1.0pt] (v2k-1) to (v2k-3);

\node[vertex] (x1) at  (0,0) {$0$};
\node[vertex] (x2) at  (0,-3) {$0$};
\node[vertex] (x3) at  (4,0) {$1$};
\node[vertex] (x4) at  (4,-3) {$1$};
\node[vertex] (x2k-2) at  (8,0) {$0$};
\node[vertex] (x2k-1) at  (8,-3) {$0$};
\node[vertex] (x2k) at  (12,0) {$1$};
\node[vertex] (x2k-3) at  (12,-3) {$1$};
\node[vertex] (x2k-4) at  (12,-1.5) {$1$};

\draw[edge,-, line width=1.0pt] (x1) to (x3);
\draw[edge,-, line width=1.0pt] (x1) to (x4);
\draw[edge,-, line width=1.0pt] (x2) to (x1);
\draw[edge,-, line width=1.0pt] (x2k-2) to (x2k-4);
\draw[edge,-, line width=1.0pt] (x2k-2) to (x2k);
\draw[edge,-, line width=1.0pt] (x2k-1) to (x2k-3);

\end{tikzpicture}

\end{center}

  \caption{ $(2,3)$-orientable labelings of four 3-edge graphs. }\label{3edge}

\end{figure}

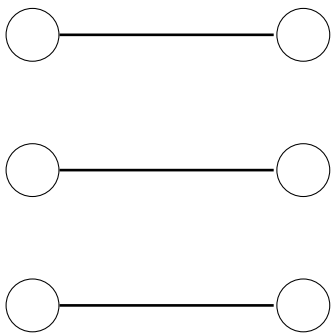
\begin{figure}[h]
\begin{center}
\begin{tikzpicture}[scale=0.9]

\tikzset{vertex/.style = {shape=circle,draw,minimum size=2em}}
\tikzset{edge/.style = {->,> = stealth',shorten >=1pt,thick}}

\node[vertex] (y1) at  (4,4) {$ $};
\node[vertex] (y2) at  (8,4) {$ $};
\node[vertex] (y3) at  (4,2) {$ $};
\node[vertex] (y4) at  (8,2) {$ $};
\node[vertex] (y5) at  (4,0) {$ $};
\node[vertex] (y6) at  (8,0) {$ $};

\draw[edge,-, line width=1.0pt] (y1) to (y2);
\draw[edge,-, line width=1.0pt] (y3) to (y4);
\draw[edge,-, line width=1.0pt] (y5) to (y6);

\end{tikzpicture}

\end{center}

  \caption{The non $(2,3)$-orientable 3-edge graph. }\label{3notedge}

\end{figure}

Note that the graph consisting of three mutually non incident arcs is $(2,3)$-cordial in the set of digraphs on seven or more  vertices but not in the set of digraphs on six vertices.  Thus the need for using nonisolated vertices  in the above definition to avoid ambiguity.  All other graphs consisting of the union of three edge graphs are $(2,3)$-orientable.  There are a total of five non isomorphic graphs of three edges.  See Figure \ref{3edge}  for the $(2,3)$-orientable ones and Figure \ref{3notedge}   for the non orientable one.

\subsubsection{Largest $(2,3)$-Orientable/Non-orientable Graphs}  

It was shown in \cite{B} that the only complete graphs that are $(2,3)$-orientable are for $n\leq 5$.  For $n \geq 6$, Santana in \cite{Sa} showed that for given $n\geq 6$ the maximum number of edges in a $(2,3)$-orientable graph is 
\begin{equation}
    \displaystyle{n\choose2} - \Big\{\displaystyle{\lceil \frac{n}{2}  \rceil\choose 2} + \displaystyle{\lfloor \frac{n}{2}  \rfloor\choose 2}\Big\} + \Bigg \lceil \frac{1}{2}\Big ( \displaystyle{n\choose2} -  \Big\{\displaystyle{\lceil \frac{n}{2}  \rceil\choose 2} + \displaystyle{\lfloor \frac{n}{2}  \rfloor\choose 2}\Big\}\Big)  \Bigg \rceil.\end{equation}  The proofs are similar to those in Section \ref{L+} above.

\section{Linear operators on graphs.}

Let $T:\Gn\to\Gn$ be a mapping such that the image of a union of graphs is the union of the images of the graphs and that the image of the edgeless (empty) graph, $\overline{K_n}$, is the edgeless graph.  Such a mapping is said to be a {\em linear operator}.  Note that $\Gn$ is a semimodule over the set $\{0,1\}$ where addition is union and  scalar multiplication is Boolean, that is for $X\in\Gn$, $0\cdot X=\overline{K_n}$ and $1\cdot X=X$.  For convenience we use juxtaposition for multiplication and we denote $\overline{K_n}$ by $O$.  So $T$ is a linear operator if for any $U,V\in\Gn$, $T(U\cup V)=T(U)\cup T(V)$ and $T(O)=O$.

Let $\X$ be a subset of $\Gn$.  Then $T$ is said to {\em preserve} the set $\X$ if whenever  $U\in\X$, $T(U)\in\X$.  We say $T$ {\em strongly preserves} $\X$ if $T$ preserves $\X$ and $T$ preserves $\Gn\setminus\X$.

\l \label{inj=bij}  Let $T:\Gn\to\Gn$ be a linear operator.  Then the following are equivalent:\begin{enumerate}\item $T$ is injective;\item $T$ is surjective;\item $T$ is bijective.\end{enumerate}\ll\d This follows from the fact that $\Gn$ is finite.\dd

Since a vertex permutation merely changes the order of the vertices (an independent labeling), we have:

\l\label{vert}  Every vertex permutation strongly preserves: \begin{enumerate}[a{\rm )}]\item  the set of $(\Z_2,+)$-cordial graphs;  \item   the set of  $(\Z_2,\times)$-cordial graphs; and \item   the set of $(2,3)$-orientable graphs.\end{enumerate}\ll

Note that a linear mapping on $\Gn$ is nonsingular if and only if the only graph mapped to the edgeless graph  is the edgeless graph.  Because this semimodule uses union for addition, nonsingularity does not imply invertibility as in vector spaces over fields.

\begin{example}
Let 
$T:{\cal G}_{10}\to {\cal G}_{10}$ and suppose that  
  $T(X)={\mathcal Pet}$ for all $X\in {\cal G}_{10}$ where ${\mathcal Pet}$ is the Petersen graph.   Then $T$  maps each 3-regular graph to a 3-regular graph ($T$ preserves 3-regular graphs).  But, $T$ also maps every non 3-regular graph to a 3-regular graph.  
\end{example}

Thus, in view of the above example, we place more restrictions on $T$ in order to reasonably characterize preservers.  This additional restriction is usually that $T$ preserves some other set or function, that $T$ is bijective, or, as below, that $T$ strongly preserves the set.

\l \label{nonsing} Let $T:\Gn\to\Gn$ be a linear operator that strongly preserves  \begin{enumerate}[a{\rm )}]\item the set of  $(\Z_2,+)$-cordial graphs, \item the set of  $(\Z_2,\times)$-cordial graphs, or \item the set of  $(2,3)$-orientable graphs.\end{enumerate} Then $T$ is nonsingular.\ll \d Suppose that $T(X)=O$ for some $X\in\Gn$.  Then for any edge graph $E\dom X$ we have that $T(E)=O$.  Let $H$ be the graph of: a) two parallel edges; b) a four cycle; or c)  three parallel edges;  one edge of which is $E$.  Then, since $T(E)=O$, $T(H)=T(H\setminus E)$, a contradiction since $H$ is not:  a) $(\Z_2,+)$-cordial; b) $(\Z_2,\times)$-cordial; or c) $(2,3)$-orientable while $H\setminus E$ is: a) $(\Z_2,+)$-cordial; b) $(\Z_2,\times)$-cordial: or c)  $(2,3)$-orientable; respectively.  Thus, in each case,  $T$ is nonsingular. \dd

\l\label{cell2cell} Let $L:\Gn\to\Gn$ be an idempotent  linear operator that strongly preserves  \begin{enumerate}[a{\rm )}]\item the set of  $(\Z_2,+)$-cordial graphs for $n\geq 4$, \item the set of  $(\Z_2,\times)$-cordial graphs for $n\geq 4$, or \item the set of  $(2,3)$-orientable graphs for $n\geq 6$.\end{enumerate}  Then the image of an edge graph is an edge graph.\ll  
\d Let $E$ be an edge graph and suppose that $L(E)\doms F$ where $F$ is an edge graph, $F\neq E$.  Say $L(E)= Q\cup F$.  Let $G$ be: a)    a $(\Z_2,+)$-cordial graph;  b) a $(\Z_2,\times)$-cordial graph; or c)   a $(2,3)$-orientable graph;
  with the maximal number of edges  such that $G\doms E$ and $G\not\doms F$.  This choice is possible since for $n\geq 4$, the complete graph $K_n$ is not  $(\Z_2,+)$-cordial  or  $(\Z_2,\times)$-cordial, and for  $n\geq 6$, $K_n$ is not  $(2,3)$-orientable.  Then $L(G)=L(G\cup E)= L(G)\cup L(E)$.  Then, $L(G)=L(G)\cup L(E)= L(G)\cup Q\cup F$, so that $L(G)=L^2(G)=L(L(G))=L(L(G)\cup L(Q)\cup L(F)) =L^2(G)\cup L(F) \cup L(Q)= L(G)\cup L(F) \cup L(Q)= L(G\cup F  \cup  Q)$, a contradiction since $G$ is:  a)  a $(\Z_2,+)$-cordial graph;  b) a $(\Z_2,\times)$-cordial graph; or c) a $(2,3)$-orientable graph,   while $G\cup F\cup Q$ cannot be: a) a  $(\Z_2,+)$-cordial graph;  b) a $(\Z_2,\times)$-cordial graph; or c) a $(2,3)$-orientable graph; respectively,    since it has too many edges. Thus $L$ maps edge graphs to edge graphs. \dd

\l\label{bij} Let $L:\Gn\to\Gn$ be an idempotent linear operator that strongly preserves  \begin{enumerate}[a{\rm )}]\item the set of  $(\Z_2,+)$-cordial graphs for $n\geq 4$, \item the set of  $(\Z_2,\times)$-cordial graphs for $n\geq 4$, or \item the set of  $(2,3)$-orientable graphs for $n\geq 6$.\end{enumerate}    Then $L$ is bijective. \ll \d  By Lemma \ref{cell2cell}, $L$ maps edge graphs to edge graphs.  Thus, $L$ is not bijective if and only if the image of two edge graphs are equal.  Suppose that $E$ and $F$ are distinct edge graphs such that $L(E)=L(F)$.  Let $G$ be: a)  a $(\Z_2,+)$-cordial graph;  b) a $(\Z_2,\times)$-cordial graph; or c) a $(2,3)$-orientable graph; with the maximal number of edges such that $G\doms E$ and $G\not\doms F$.   Then $L(G)=L(G\cup E)=L(G)\cup L(E)=L(G)\cup L(F)=L(G\cup F)$ a contradiction since 
 $G$ is:  a)  a $(\Z_2,+)$-cordial graph;  b) a $(\Z_2,\times)$-cordial graph; or c) a $(2,3)$-orientable graph,    while $G\cup F$ cannot be: a) a  $(\Z_2,+)$-cordial graph;  b) a $(\Z_2,\times)$-cordial graph; or c) a $(2,3)$-orientable graph; respectively,  since it has too many edges.. Thus $L$ is bijective on the set of edge graphs and since $L$ is linear $L$ is bijective on $\Gn$.\dd

An interesting fact about idempotent maps on finite sets is that they are bijective if and only if they are the identity:

\l  Let ${\cal K}$ be a finite set and let $\Phi:{\cal K}\to{\cal K}$.  If $\Phi$ is bijective and idempotent then $\Phi$ is the identity map.\ll

\d Assume that $\Phi$ is bijective.   Let $A\in{\cal K}$ and suppose that $\Phi(A)=B$.  Then,
 $\Phi(B)\,\,\,=\Phi(\Phi(A)) = \Phi^2(A) =\Phi(A)$  (since $\Phi$ is idempotent).  
So that $\Phi(A)=\Phi(B)$.  Therefore $A=B$ since $\Phi$ is bijective and hence, $\Phi(A)=A$ for all $A\in{\cal K}$. \dd

Note that because union is the addition in $\Gn$, any bijective map on $\Gn$ maps edge graphs to edge graphs and preserves $|E(G)|$, the number of edges in a graph 

The next theorem uses the following lemma  from \cite{BP}:
\l\label{verper}{\rm \cite[Lemma 2.2]{BP}} If $T:\Gn\to\Gn$ is bijective, preserves $|E(G)|$, and maps 2-stars to 2-stars then $T$ is a vertex permutation.\ll

\t  Let $T:\Gn\to\Gn$ be a linear operator.  Then, $T$  strongly preserves  \begin{enumerate}[a{\rm )}]\item the set of  $(\Z_2,+)$-cordial graphs for $n\geq 4$, \item the set of  $(\Z_2,\times)$-cordial graphs for $n\geq 5$, or \item the set of  $(2,3)$-orientable graphs for $n\geq 6$\end{enumerate} \vskip -8pt  \quad  if and only if \\  $T$ is a vertex permutation.\tt \d  By Lemma \ref{vert},  every vertex permutation strongly preserves the set of  $(\Z_2,+)$-cordial graphs, the set of  $(\Z_2,\times)$-cordial graphs and  the set of $(2,3)$-orientable graphs.  So we now assume that $T$   strongly preserves: a)  the set of  $(\Z_2,+)$-cordial graphs; b)the set of  $(\Z_2,\times)$-cordial graphs, or c) the set of  $(2,3)$-orientable graphs.  Let $L=T^d$ where $d$ is chosen so that $L$ is idempotent.  By Lemma \ref{bij}, $L$ is bijective and hence the identity mapping.    Now, suppose that $T(X)=T(Y)$.  Then $T^{d-1}(T(X))=T^{d-1}(T(Y))$ so that $L(X)=T^d(X)= T^{d-1}(T(X))=T^{d-1}(T(Y))=T^d(Y)=L(Y)$.  Since $L$ is the identity, we have $X=Y$.  That is $T$ is injective and hence, by Lemma \ref{inj=bij}, $T$ is bijective.  Since addition is union, the image of any edge graph is a edge graph.  Thus, for any graph, $|E(G)|=|E(T(G))|$.  That is, $T$ preserves $|E(G)|$.

{\large\bf Case 1.} $T$ preserves the set of $(\Z_2,+)$-cordial; graphs.  Let $E$ and $F$ be incident edges, so that $E\cup F$ is a 2-star (2-path) and suppose that $T(E\cup F)$ is a pair of non incident (parallel) edges.  Then $E\cup F$ is $(\Z_2,+)$-cordial while its image is not, a contradiction.   That is $T$ maps 2-stars to 2-stars.  By Lemma \ref{verper}, $T$ is a vertex permutation.

{\large\bf Case 2.} $T$ preserves the set of $(\Z_2,\times)$-cordial; graphs.    Let $E$ and $F$ be incident edges, so that $E\cup F$ is a 2-star (2-path) and suppose that $T(E\cup F)$ is a pair of non incident (parallel) edges.  Let $K$ be the edge joining the vertices of $E\cup F$ that forms a triangle.  Let $G=E\cup F\cup K$.   Then there are $3(n-3)$ edges that are adjacent to G.  There are only three possible images of $G$ since $T$ is bijective and  the image of any edge graph is an edge graph. $T(G)$ is either a 2-star with a non-incident edge, a three path or three mutually parallel edges (a $3K_2$ ).  If $T(G)$ is a 2-star with a non-incident edge, there is only one edge  added to $T(G)$ that results in a non $(\Z_2,\times)$-cordial graph.  If $T(G)$ is a $3K_2$, no edge added to $T(G)$ results in a non $(\Z_2,\times)$-cordial graph, and if $T(G)$ is a three path, there are exactly three edges which, when any one is added to $T(G)$ would result in a  non $(\Z_2,\times)$-cordial graph.  Since $G$ together with any one of the $3(n-3)$ edges adjacent to $G$ is a non $(\Z_2,\times)$-cordial graph, so must its image be.  Since we have shown that this is impossible, it follows that $T$ maps 2-stars to 2-stars. By Lemma \ref{verper}, $T$ is a vertex permutation.

{\large\bf Case 3.} $T$ preserves the set of $(2,3)$-orientable graphs.  
Let $E$ and $F$ be incident edges, so that $E\cup F$ is a 2-star (2-path) and suppose that $T(E\cup F)$ is a pair of non incident (parallel) edges.  Let $H$ be an edge graph such that the image of $H$ is an edge graph parallel to both $T(E)$ and $T(F)$.  Then, $E\cup F\cup H$ is $(2,3)$-orientable while $T(E\cup F\cup H)$, being three parallel edges, is not $(2,3)$-orientable, a contradiction.  That is $T$ maps 2-stars to 2-stars.  By Lemma \ref{verper}, $T$ is a vertex permutation.  \dd

\end{document}